\documentclass[10pt,twoside]{article}
\usepackage{Latex-document}

\newtheorem{theorem}{Theorem}[section]

\newtheorem{definition}[theorem]{Definition}

\def\Section#1{\section{\hskip -1em . \hskip 0.8em #1}}
\renewcommand{\thesection}{\arabic{section}}
\def\addsec{\setzero \setcounter{theorem}{0}}

\title{\bf Some Geometric PDEs Related to \vskip -2mm Hydrodynamics and Electrodynamics\vskip -2mm \vskip 6mm}
\author{Yann Brenier\vspace*{-0.5cm}\thanks{CNRS, Universit\'e de Nice-Sophia-Antipolis, France
(on leave of absence from Universit\'e Paris 6) and Institut Universitaire de France. E-mail:
brenier@math.unice.fr}}
\date{\vspace{-8mm}}

\begin{document}

\maketitle

\thispagestyle{first} \setcounter{page}{761}

\begin{abstract}

\vskip 3mm

We discuss several geometric PDEs and their relationship with
Hydrodynamics and classical Electrodynamics. We start from the
Euler equations of ideal incompressible fluids that, geometrically
speaking, describe geodesics on groups of measure preserving maps
with respect to the $L^2$ metric. Then, we introduce a geometric
approximation of the Euler equation, which involves the
Monge-Amp\`ere equation and the Monge-Kantorovich optimal
transportation theory. This equation can be interpreted as a fully
nonlinear correction of the Vlasov-Poisson system that describes
the motion of electrons in a uniform neutralizing background
through Coulomb interactions. Finally we briefly discuss an
equation for generalized extremal surfaces in the 5 dimensional
Minkowski space, related to the Born-Infeld equations, from which
the Vlasov-Maxwell system of classical Electrodynamics can be
formally derived.

\vskip 4.5mm


\end{abstract}

\vskip 12mm

\Section{The Euler equations of incompressible fluids} \addsec

\vskip-5mm \hspace{5mm}

The motion of
an incompressible fluid moving in a compact domain $D$ of the Euclidean
space ${{\bf R}}^d$ can be mathematically defined as a trajectory
 $t\rightarrow g(t)$ on the
set, subsequently denoted by $G(D)$,
of all diffeomorphisms of $D$ with unit jacobian determinant.
This space can be embedded in the set $S(D)$
of all Borel maps $h$ from $D$ into itself,
not necessarily one-to-one, such that
\[\int_D\phi(h(x))dx=\int_D\phi(x)dx\]
for all $\phi\in C(D)$, where $dx$ denotes the Lebesgue measure,
normalized so that the measure of $D$ is 1. For the composition
rule, $G(D)$ is a group (the identity map $I$ being the unity of
the group), meanwhile $S(D)$ is a semi-group. Both $G(D)$ and
$S(D)$ are naturally embedded in the Hilbert space $H=L^2(D,{{\bf
R}}^d)$ of all square integrable mapping from $D$ into ${{\bf
R}}^d$ and, therefore, inherit from $H$ a formal Riemannian
structure. The equations of geodesics on $G(D)$ turn out to be
exactly \cite{AK} the equations of incompressible inviscid fluids
introduced by Euler near 1750 \cite{Eu}. The Euler equations play
a fundamental role in Fluid Mechanics (for geophysical flow
modelling in particular) and their global well-posedness is one of
the most challenging problems in the field of nonlinear PDEs.
Their mathematical importance is confirmed by the recent
publication of several books by Arnold-Khesin \cite{AK}, Chemin
\cite{Ch}, P.-L. Lions \cite{Li}, Marchioro-Pulvirenti \cite{MP},
as well as by Majda's lecture in the Kyoto ICM \cite{Ma}.

From a geometric point of view (different from the usual PDE
setting which consists in solving the Euler equations with
prescribed initial conditions), it is natural to look for
minimizing geodesics between the identity map and prescribed
measure preserving maps. More precisely~:

\begin{definition}
Given $h\in G(D)$, find a curve
$t\in [0,1]\rightarrow g(t)\in G(D)$ satisfying $g(0)=I$, $g(1)=h$,
that minimizes
\[A_D(g)=\frac{1}{2}\int_0^1||g'(t)||_{L^2}^2 dt
=\frac{1}{2}\int_0^1\int_D|\partial_t g(t,x)|^2 dxdt.\]
\end{definition}

The infimum is nothing but $\frac{1}{2}\delta_D^2(I,h)$, where
$\delta_D$ denotes the geodesic distance on $G(D)$, and any smooth
minimizer $g$ must be a smooth solution of  the Euler equations
(written in ``Lagrangian coordinates'')
\[g''\circ g^{-1}=-\nabla p,\]
where $p=p(t,x)\in\bf{R}$ is the pressure field and
$\nabla p=(\partial_{x_1}p,...,\partial_{x_d}p)$.
The minimization problem will be subsequently called
``Shortest Path Problem'' (SPP).
\\
The basic local existence and uniqueness theorem for the SPP
is due to Ebin and Marsden \cite{EM}.
If $h$ and $I$ are sufficiently close in a sufficiently high order Sobolev
norm, then there is a unique shortest path. In the large, uniqueness can
fail for the SPP. For example, in the case when $D$ is the unit disk,
$h(z)=-z$, the SPP has two solutions $g(t,z)=ze^{+i\pi t}$
and $g(t,z)=ze^{-i\pi t}$, where complex notations are used.
\\
In 1985, A. Shnirelman \cite{Sh} found, in the case $D=[0,1]^3$,
a class of data for which the
SPP cannot have a (classical) solution.
These data are those of form
\[h(x_1,x_2,x_3)=(H(x_1,x_2),x_3),\]
where $H$ is an area preserving mapping of the unit square, i.e. an element
of $G([0,1]^2)$, for which
\[\delta_{[0,1]^3}(I,h)<\delta_{[0,1]^2}(I,H)<+\infty.\]
(This means that, although $h$ is really a two dimensional map,
genuinely 3D motions perform better to reach $h$ from $I$
than purely 2D motions.)
\\
Shnirelman
also proved \cite{Sh}, \cite{Sh2},
that $S([0,1]^d)$ is the right completion of $G([0,1]^d)$
for the geodesic distance $\delta$, for all dimension $d\ge 3$.
(Notice that $S([0,1]^d)$ is the $L^2$ completion of $G([0,1]^d)$
for all $d\ge 2$ \cite{Ne}. So, the case $d=2$ is very peculiar.)
\\
In such situations, a complete existence and uniqueness result for
the SPP was obtained in \cite{Br2},
$provided$ the pressure field is considered as the right unknown
and not the path $t\rightarrow g(t)$ itself.

\begin{theorem}
\label{detailF.1}

Let $h\in S([0,1]^3)$ of form $h(x_1,x_2,x_3)=(H(x_1,x_2),x_3)$
with $H\in S([0,1]^2)$.
Then there is a unique vector-valued measure
$\nabla p(t,x_1,x_2)$ such that, for each
sequence of curves $t\in [0,1]\rightarrow g_n(t)\in G([0,1]^3)$
labelled by $n\in \bf{N}$ and satisfying
\[
A_{[0,1]^3}(g_n)\rightarrow \frac{1}{2}\delta_{[0,1]^3}^2(I,h),\;\;\;
||g_n(1)-h||_{L^2([0,1]^3}\rightarrow 0,
\]
as $n\rightarrow \infty$, then (in the distributional sense)
\[
g_n''\circ g_n^{-1}\rightarrow -\nabla p.
\]
\end{theorem}
In other words, the acceleration field of $all$ minimizing sequences
converge to $-\nabla p$ which uniquely depends on data $h$.
The proof relies on an appropriate concept of generalized
solutions (related to ``Young's measures'' ~\cite{Yo}, \cite{Ta},
\cite{DM}, \cite{She})
that describe the oscillatory behaviour of
the $(g_n)$ as $n\rightarrow +\infty$ and reduces the SPP
to a convex minimization problem. (See \cite{Br2} for
more details.)
More precisely, the associated measures
\[
c_n(t,x,a)=\delta(x-g_n(t,a)),\;\;\;
m_n(t,x,a)=\partial_t g_n(t,a)\delta(x-g_n(t,a)),
\]
have cluster points $(c,m)$ that have the following properties~:
\\
1) $m$ is absolutely continuous with respect to $c$ and its vector-valued
density $v(t,x,a)$ is $c-$ square integrable;
\\
2) $c$ and $v$ do not depend on $x_3$ and $v_3=0$,
\\
3) $c$ and $v$ solve
\begin{equation}
\label{hydrostatique}
\partial_t c+\nabla_x.(cv)=0,\;\;\;
\partial_t (cv)+\nabla_x.(cv\otimes v)+c\nabla_x p=0,
\end{equation}
where the product $c\nabla_x p$ has to be properly defined (in a
way related to the work of Zheng and Majda \cite{ZM}). Equations
(\ref{hydrostatique}) are obtained as the optimality equations of
the convexified minimization problem. Therefore, it is a priori
unclear they have any physical meaning as evolution equations.
However, they correspond, up to a change of unknown, to the
hydrostatic limit of the Euler equations, obtained from the Euler
equations by neglecting the vertical acceleration term, namely:
\[Kg''\circ g^{-1}=-\nabla p,\]
where $K$ is the singular diagonal matrix $(1,1,0)$. These hydrostatic
(or ``shear flow'')
equations are widely used for atmosphere and ocean circulation modelling,
as the building block of the so-called ``primitive equations''.
However, they are more singular than the Euler equations and their
mathematical analysis is very limited, as discussed in \cite{Li}.
Conditional well posedness and derivation from the
Euler equations have been established in \cite{Br3} and \cite{Gr}.

\subsubsection*{Remarks}

\vskip -5mm \hspace{5mm}

 An intriguing question is whether or not the
uniqueness of $\nabla p$ can be proved by more classical tools
even in the case when $H\in G([0,1]^2)$ can be connected to the
identity map by a classical shortest path on $G([0,1]^2)$.

Since $S([0,1]^3)$ is the right completion of $G([0,1]^3)$ with
respect to the geodesic distance, one could expect the $SPP$ to
have a solution in $S([0,1]^3)$ for all data $h$. This is not true.
An example of such a data is $h(x_1,x_2,x_3)=(1-x_1,x_2,x_3)$.
Only generalized flows, as discussed in \cite{Br2}, can describe
shortest paths in full generality.

\subsection*{Example of generalized solutions}

\vskip -5mm \hspace{5mm}

Explicit examples of non trivial generalized shortest paths can be
computed either numerically or exactly.  Let us just quote a
typical example, when $D$ is the cylinder
$\{(z,s)=(x_1,x_2,s),\;\;|z|\le 1,\;\;0\le s\le 1\}$ and
$h(z,s)=(-z,s)$. Then, the classical SPP has two distinct
solutions $g_+(t,z,s)=(e^{i\pi t}z,s)$ and $g_-(t,z,s)=(e^{-i\pi
t}z,s)$, with the same pressure field $p=\pi^2|z|^2/2$, where
complex notations are used on the disk $|z|\le 1$. (Notice that
there is no motion along the vertical axis $s$.) Trivial
generalized solutions are obtained by mixing these two solutions.
However, a non trivial generalized solution exists and can be
described as follows. For each fluid particle initially located at
$(z,s)$, the elevation $s$ stays unchanged and the initial
horizontal position $z$ splits up along a circle of radius
$(1-|z|^2)^{1/2}\sin(\pi t)$, with center $z\cos(\pi t)$, that
moves across the unit disk and shrinks down to the point $-z$ as
$t=1$. In addition, each particle is accelerated by the pressure
field $p=\pi^2|z|^2/2$, as expected from the theory.

\Section{Polar factorization of maps and the Monge-\\Amp\`ere equation}
\addsec

\vskip-5mm \hspace{5mm}

A way to define approximate geodesics on $G=G(D)$ is to
introduce a penalty parameter $\epsilon>0$ and to consider
the formal (hamiltonian) dynamical system in the Hilbert space
$H=L^2(D,{{\bf R}}^d)$
\begin{equation}
\label{hamiltonien0}
\epsilon^2\frac{d^2}{dt^2}M+\frac{\delta}{\delta M}\left( \frac{d_{_H}^2(M,G)}
{2}\right)
=0,
\end{equation}
where the unknown $M$ is a curve $t\rightarrow M(t)\in H$,
$\delta/\delta M$ denotes the gradient operator in $H$, and
\begin{equation}
\label{potentiel}
d_{_H}(M,G)
=\inf_{g\in G}||M-g||_{_H}
\end{equation}
is the distance in $H$ between $M$ and $G$,
where $||.||_{_H}$ is the Hilbert norm of $H$.
This approach is related
to Ebin's slightly compressible flow theory \cite{Eb}, and is
a natural extension of the theory of
constrained finite dimensional mechanical systems \cite{RU}, \cite{AK}.
Notice that the approximate geodesic equation is sensitive only to
the $L^2$ closure of $G(D)$, which is, in the case
$D=[0,1]^d$, $d\ge 2$, the entire semi-group $S(D)$ \cite{Ne}.
As the penalty parameter $\epsilon$ goes to zero,
we expect that for appropriate initial data,
typically for $M(t=0)=I$ and $(d/dt)M(t=0)=v_0$, where $v_0$ is
a smooth divergence free vector field on $D$ tangent to the boundary,
the time dependent map $M$ converges to a geodesic curve on $G$.
Because of the classical properties of the distance function in a
Hilbert space, for each point $M\in H$ for which there exists
a unique closest point $\pi_{S}(M)$ on $S(D)$, we have
\begin{equation}
\frac{\delta}{\delta M}( \frac{d_{_H}^2(M,G)}{2})=M-\pi_S (M).
\end{equation}
Thus, we can formally write
the approximate geodesic equation (\ref{hamiltonien0})
\begin{equation}
\label{hamiltonien}
\epsilon^2\frac{d^2}{dt^2}M+M-\pi_S (M)=0.
\end{equation}
Therefore, it is natural to address the following variational problem,
that we call the Closest Point Problem (CPP)

\begin{definition}
Given $M\in L^2(D,{{\bf R}}^d)$,
find $h\in S(D)$ that minimizes
\[\frac{1}{2}\int_D|M(x)-h(x)|^2 dx.\]
\end{definition}

The solution of the CPP is given by the Polar Factorization theorem
for maps \cite{Br1}

\begin{theorem}
Let $M:D\rightarrow {{\bf R}}^d$ be an $L^2$ map
such that the probability measure
\[
\rho_M(x)=\int_D \delta(x-M(a))da
\]
is a Lebesgue integrable function on $D$. Then, there exists a unique
closest point $\pi_S(M)$ on $S(D)$ and there is a Lipschitz
convex function $\Phi$ on ${{\bf R}}^d$ such that
\[
\pi_S(M)(a)=(\nabla\Phi)(M(a)),\;\;a.e.\;a \in D.
\]
In addition, $\Phi$ is a weak solution, in a suitable sense,
of the Monge-Amp\`ere equation
\[
\det(\partial_{xx}\Phi(x))=\rho_M(x).
\]
\end{theorem}

Thus, the Monge-Amp\`ere equation \cite{Ca}, which is usually considered
as a non variational geometric PDE related to the concept of
Gaussian curvature, also is the optimality equation of a variational
problem closely linked to the Euler equations of incompressible
inviscid fluids.
In addition, the Polar Factorization theorem can be seen as a nonlinear
version of the Helmholtz-Hodge decomposition theorem for vector fields
which asserts that any $L^2$ vector field on $D$ can be written in a unique
way as the (orthogonal) sum of the gradient of
a scalar field and a divergence free field tangent to $\partial D$.
Shortly after \cite{Br1},
Caffarelli \cite{Ca}
established several regularity results for the Polar Factorization.
For example, provided $D$ is smooth and strictly convex,
any smooth orientation preserving diffeomorphism $M$ of $D$ has a unique
Polar Factorization with smooth factors, and
$\pi_S(M)$ belongs to $G(D)$. More recently, McCann \cite{Mc}
generalized the Polar Factorization theorem
when $D$ is a compact Riemannian manifolds.

\Section{Optimal Transportation Theory} \addsec

\vskip-5mm \hspace{5mm}

In \cite{Br1}, the solution of the CPP problem
is based on the Optimal Transportation Theory (OTT).
The OTT was introduced by Monge in 1781 \cite{Mo}
to solve an engineering problem
and renewed by Kantorovich near 1940 \cite{Ka} in the framework of
Linear Programing and Probability Theory \cite{RR}. In modern words, this
amounts to look for a probability measure $\mu$ on a given product
measure space $A\times B$, with prescribed projections on $A$ and $B$,
that minimizes
\[
\int_{A\times B}c(x,y)d\mu(x,y),
\]
where the ``cost function'' $c\ge 0$ is given on $A\times B$.
The CPP roughly corresponds to the case when $A=B=D$,
$c(x,y)=|M(x)-y|^2$ and each projection of $\mu$ is the (normalized)
Lebesgue measure on $D$.
The connexion established in \cite{Br1}
between the OTT and the Monge-Amp\`ere equation, enhanced
by Caffarelli's regularity theory \cite{Ca},
introduced OTT as an active field of research in nonlinear PDEs.
Let us first quote the work of Evans-Gangbo \cite{Ev}
to solve the original Monge problem with PDE techniques, related to the
Eikonal equations,
and the recent contributions of Ambrosio, Caffarelli, Feldman, McCann,
Trudinger, Wang. (A first attempt was made by Sudakov \cite{Su} with
purely probabilistic tools.)
Let us next point out the importance of OTT for modelling
purposes in Applied Mathematics. First of all, it is fair to say that
the OTT and the Monge-Amp\`ere equation were
already key ingredients in Cullen and Purser's theory of
semi-geostrophic atmospheric flows, which goes back to the early 1980s and
preceded our Polar Factorization theorem (see references in \cite{CNP}).
Next, Jordan-Kinderlehrer-Otto \cite{JKO}, using OTT,
established that the heat
equation can be seen as a gradient flow for Boltzmann's entropy functional.
More systematically, Otto \cite{Ot} showed how the OTT confers
a natural Riemannian structure to sets of Probability measures
and recognized a large class of dissipative PDEs as
gradient flows of various functionals for such Riemannian structures.
Examples of such PDEs are porous media equations,
lubrication equations, granular flow equations, etc...
Let us also mention that OTT has became a powerful tool in Calculus
of Variations (through McCann's concept of displacement convexity
\cite{Mc})
and Functional Analysis, where all kind of functional inequalities
(Minkowski, Brascamp-Lieb, Log Sobolev,
Bacry-Emery, etc,...) can be established through OTT arguments,
as shown, in particular,
by Barthe \cite{Ba}, McCann \cite{Mc}, Otto, Villani \cite{OV}.
Let us finally mention that \cite{BB} has provided for the OTT
a formulation different from the Monge-Kantorovich one,
by introducing an interpolation variable
(which was already present in McCann's concept of displacement
convexity). This point of view
is useful for both numerical \cite{BB}
and theoretical purposes, in particular,  by allowing non trivial
generalizations of the OTT related to section \ref{courant}.

\Section{Approximate geodesics and Electrodynamics} \addsec

\vskip-5mm \hspace{5mm}

Let us go back to the approximate geodesic equation
(\ref{hamiltonien}) that can be
(formally) written, thanks to the Polar Factorization theorem,
\begin{equation}
\label{hamiltonien bis}
\partial_{tt}M(t,a)+(\nabla\phi)(t,M(t,a))=0,
\;\;\;\det(I-\epsilon^2 \partial_{xx}\phi(t,x))=\rho_M(t,x)
\end{equation}
(where $\phi(t,x)$ stands for $\epsilon^{-2}(|x|^2/2-\Phi(t,x))$).
A formal expansion about $\epsilon=0$ leads, as expected, to
the Euler equation (written in Lagrangian coordinates)
at the zero order and,
at the next order (and exactly as $d=1$), to
\begin{equation}
\label{poisson}
\partial_{tt}M(t,a)+(\nabla\phi)(t,M(t,a))=0,
\;\;\;\epsilon^2 \Delta\phi(t,x)=1-\rho_M(t,x),
\end{equation}
which can be equivalently written as
\begin{equation}
\label{vlasov poisson}
\partial_{t}f+\xi.\nabla_x f-\nabla_x\phi.\nabla_\xi f=0,
\;\;\;\epsilon^2 \Delta\phi=1-\int fd\xi
\end{equation}
by introducing the ``phase density''
\[
f(t,x,\xi)=\int_D \delta(x-M(t,a))\delta(\xi-\partial_t M(t,a))da.
\]
This system is nothing but the Vlasov-Poisson system
that describes the classical non-relativistic motion of a
continuum
of electrons around a homogeneous neutralizing background of ions
through Coulomb interactions.

So, the approximate geodesic equation, which can be written
as a ``Vlasov-Monge-Amp\`ere'' (VMA) system,
\begin{equation}
\label{VMA}
\partial_{t}f+\xi.\nabla_x f-\nabla_x\phi.\nabla_\xi f=0,
\;\;\;\det(I-\epsilon^2 \partial_{xx}\phi)=\int fd\xi
\end{equation}
can be interpreted as a (fully nonlinear) correction of the Vlasov-Poisson
system for small values of $\epsilon$.
Recently, Loeper \cite{Lo}
has shown that the VMA system has local smooth solutions
and global weak solutions. Loeper has also proved that
the Euler equations and the Vlasov-Poisson system correctly
describe the asymptotic
behaviour of the VMA system as $\epsilon\rightarrow 0$. The asymptotic
analysis is based on the so-called
modulated energy method already
used in \cite{Br5} to derive the Euler equations
from the Vlasov-Poisson system.

Notice that, thanks to the substitution of the Monge-Amp\`ere equation
(a fully non-linear elliptic PDE) for the classical Poisson equation,
the ``electric'' field $\nabla\phi(t,x)$ is pointwise bounded by
the diameter of $D$ divided by $\epsilon^2$, independently on the initial
conditions. In particular, point charges do not create unbounded
force fields as in classical Electrodynamics.

\Section{A caricature of Coulomb interaction} \addsec

\vskip-5mm \hspace{5mm}

The approximate geodesic equation (\ref{hamiltonien}) can be
easily discretized in space by substituting i) for $D$ a discrete
set of $N$ ``grid'' points equally spaced in $D$, say
$A_1,...,A_N$,  ii) for $H$ the euclidean space ${\bf R}^{dN}$,
iii) for $G$ the discrete set of all sequences
$(A_{\sigma_1},...,A_{\sigma_N})\in {\bf R}^{dN}$ generated by
permutations $\sigma$ of the first $N$ integers, while keeping
unchanged equation (\ref{hamiltonien}). (Note that such a
discretization using permutations cannot be so easily defined for
the Euler equations, which formally correspond to the limit case
$\epsilon=0$.) Then $M(t)=(M_1(t),...,M_N(t))$ can be interpreted
as a set of $N$ harmonic oscillators
\begin{equation}
\label{matchon}
\epsilon^2 \frac{d^2}{dt^2}M_\alpha+M_\alpha-A_{\sigma_\alpha(t)}=0,
\end{equation}
where the time dependent permutation $\sigma(t)$ is subject to minimize,
at all time $t$, the total potential energy
\begin{equation}
\label{matchon bis}
\sum_{\alpha=1}^N |M_\alpha(t)-A_{\sigma_\alpha}|^2.
\end{equation}
This system can be seen as a collection of $N$ springs
linking each particle $M_\alpha$ to one of the fixed particle $A_\beta$
according to a dynamical pairing $\beta=\sigma_\alpha(t)$
maintaining the bulk potential energy at the lowest level.
There is some ambiguity in the definition of this formal hamiltonian
system for which the hamiltonian is given by
\begin{equation}
\label{matchon ter} \frac{1}{2}\sum_{\alpha=1}^N
|\frac{dM_\alpha}{dt}|^2
+\inf_{\sigma}\frac{1}{2\epsilon^2}\sum_{\alpha=1}^N
|M_\alpha-A_{\sigma_\alpha}|^2.
\end{equation}
In particular, $\sigma(t)$ is not uniquely defined at each time $t$
for which several particles have the same position.
However, the potential is the sum of a quadratic and a Lipschitz concave
functions of $M$. So its gradient has linear growth at infinity and its
second order partial derivatives are locally bounded measures.
This is enough, according to recent results by Lions and Bouchut
\cite{Bo}, \cite{Li2}, to ensure
that unique global solutions are well defined for
Lebesgue almost every initial data $M_\alpha(0)$, $\frac{d}{dt}M_\alpha(0)$,
$\alpha=1,...,N$.
As expected, the limit $N\rightarrow +\infty$, $\epsilon\rightarrow 0$
(provided $N$ goes fast enough to $+\infty$), leads to the Euler
equation, as proven in \cite{Br4}.
From the electrostatic point of view,
the dynamical system describes a nonlinearly cutoff Coulomb
interaction between $N$ electrons (with positions $M_\alpha$)
and a background of $N$ motionless ions (with fixed positions $A_\alpha$).

\Section{Generalized extremal surface equations and Electrodynamics}
\label{courant} \addsec

\vskip-5mm \hspace{5mm}

As seen above, the approximate geodesic equation (\ref{hamiltonien
bis})---which has been introduced as a natural geometrical
approximation to the Euler equations---turns out to be a model for
electrostatic interaction with a non-linearly cutoff Coulomb
potential. This feature is somewhat reminiscent of the Born-Infeld
non-linear theory of the electromagnetic field \cite{BI} (see also
\cite{BDLL}, \cite{GZ}...). Therefore, one may try to design from
similar geometric ideas a non-linearly cutoff theory for classical
Electrodynamics. An attempt is made in \cite{Br6}. Instead of
considering springs linking two particles of opposite charges we
rather consider (with a more space-time oriented point of view)
surfaces $(t,s)\rightarrow X(t,s)$ spanning curves $t\rightarrow
X_-(t)$ and $t\rightarrow X_+(t)$ followed by two particles of
opposite charge, so that $X(s=-1,t)=X_-(t)$ and $X(s=1,t)=X_+(t)$,
$s\in [-1,1]$ standing for the ``interpolation'' parameter between
the two trajectories. Just by prescribing $(t,s)\rightarrow
(t,s,X(t,s))$ to be an extremal surface in the the 5 dimensional
Minkowski space $(t,s,x_1,x_2,x_3)$ (with signature $(-++++)$), we
get the building block of the model. In other words, the
individual Action of each surface is
\begin{equation}
\int\sqrt{1+|\partial_s X|^2-|\partial_t X|^2-
|\partial_s X\times \partial_t X|^2}dtds,
\end{equation}
(which is basically the Nambu-Goto Action of classical string theory).
Next, we associate with $X$ a ``generalized surface'' $(\rho,J,E,B)$
(or more precisely a ``cartesian current'' in the sense of \cite{GMS})
defined by
\begin{equation}
\rho(t,s,x)=\delta(x-X(t,s)),\;\;\; J(t,s,x)=\partial_t
X(t,s)\delta(x-X(t,s)),
\end{equation}
\begin {equation}
E(t,s,x)=\partial_s X(t,s)\delta(x-X(t,s)),
\end{equation}
\begin{equation}
B(t,s,x)=\partial_s X(t,s)\times \partial_t X(t,s)\delta(x-X(t,s))
\end{equation}
and subject to compatibility conditions
\begin{equation}
\label{compatibilite}
\partial_s\rho+\nabla.E=0,\;\;\partial_t\rho+\nabla.J=0,\;\;
\partial_t E-\partial_s J-\nabla\times B=0.
\end{equation}
In terms of $(\rho,J,E,B)$ the Action of $X$ can be written as
\begin{equation}
K(\rho,J,E,B)=\int \sqrt{\rho^2-J^2+E^2-B^2}.
\end{equation}
Varying this Action under constraint (\ref{compatibilite}) leads to
a system of evolution equations for $(\rho,J,E,B)$
(see \cite{Br6} for an explicit form), that we can
call ``generalized extremal surface equations'' (GESE). They enjoy
(at least in the simplest cases when the solutions depend on one or
two space variables) many interesting properties :
hyperbolicity, linear degeneracy of all fields \cite{BDLL}, symmetries
between $t$ and $s$, $J$ and $E$ etc...
From the GESE, we can derive through various (formal!) limiting
process 1) the Born-Infeld and the Maxwell equations,
as $(\rho,J)$ are prescribed at $s=-1$ and $s=+1$ (in which case there
is no coupling between charged particles and the electromagnetic field),
2) the Vlasov-Born-Infeld and the Vlasov-Maxwell equations as
$(E,B)=0$ is prescribed at $s=-1$ and $s=+1$ (which corresponds to a free
boundary condition $\partial_s X=0$ for an individual surface and yields
a full coupling between charged particles and the electromagnetic field).
In spite of the possible physical irrelevance of the GESE,
their mathematical analysis (global existence, uniqueness, etc...),
and the rigorous derivation from them of classical models, such as
the Vlasov-Maxwell equations, are, in our opinion,
challenging problems in the field of non-linear PDEs.

\bibliographystyle{amsalpha}

\label{lastpage}

\end{document}